\documentclass[11pt]{amsart}
\usepackage[T1]{fontenc}
\usepackage{lmodern,amsmath,amssymb,amsthm,needspace}
\usepackage[a4paper,margin=28mm]{geometry}
\usepackage[hidelinks]{hyperref}
\newtheorem{theorem}{Theorem}[section]
\newtheorem{lemma}[theorem]{Lemma}
\newtheorem{proposition}[theorem]{Proposition}
\newtheorem{corollary}[theorem]{Corollary}
\newtheorem{conjecture}[theorem]{Conjecture}
\theoremstyle{definition}

\newtheorem{question}[theorem]{Question}
\theoremstyle{remark}
\newtheorem{remark}[theorem]{Remark}
\newcommand{\PF}{\operatorname{PF}}
\newcommand{\Tor}{\operatorname{Tor}}
\newcommand{\Hom}{\operatorname{Hom}}
\newcommand{\hml}{\operatorname{H}}
\newcommand{\im}{\operatorname{Im}}
\newcommand{\Ker}{\operatorname{Ker}}
\newcommand{\Ann}{\operatorname{Ann}}
\newcommand{\Ap}{\operatorname{Ap}}
\newcommand{\type}{\operatorname{type}}
\newcommand{\m}{\mathfrak m}
\newcommand{\n}{\mathfrak n}
\newcommand{\fkc}{\mathfrak c}
\newcommand{\N}{\mathbb N}
\newcommand{\detid}{\operatorname{I}}

\title[Arithmetic pseudo-Frobenius numbers and Koszul homology]{Arithmetic pseudo-Frobenius numbers and Koszul homology in numerical semigroup rings}
\author{Toshinori Kobayashi}
\address[Toshinori Kobayashi]{Department of Mathematics, School of Science and Technology, Meiji University,
1-1-1 Higashi-mita, Tama-ku, Kawasaki 214-8571, Japan}
\email{tkobayashi@meiji.ac.jp}
\author{Naoyuki Matsuoka}
\address[Naoyuki Matsuoka]{Department of Mathematics, School of Science and Technology,
Meiji University, 1-1-1 Higashi-mita, Tama-ku, Kawasaki 214-8571, Japan}
\email{naomatsu@meiji.ac.jp}
\subjclass[2020]{Primary 13C40, Secondary 13F65, 13D02}
\keywords{Numerical semigroup ring, pseudo-Frobenius number, determinantal ideal, Koszul homology, Golod ring}
\thanks{The authors were partially supported by JSPS KAKENHI Grant Numbers JP25K17240 (Toshinori Kobayashi) and JP25K06941 (Naoyuki Matsuoka).}
\thanks{The authors used OpenAI's ChatGPT to discuss proofs and to assist with drafting and revising the manuscript. The authors checked the mathematical arguments and take full responsibility for the results and the exposition.}

\begin{document}
\begin{abstract}
D. T. Cuong, D. V. Kien, N. Matsuoka, and H. L. Truong conjectured that, if a numerical semigroup $H$ has $n$ minimal generators and its pseudo-Frobenius numbers form an arithmetic progression of length $n-1$, then the defining ideal of $k[H]$ is generated by the $2\times2$ minors of a $2\times n$ matrix whose entries are powers of the variables. We prove this conjecture when $k[[H]]$ is Golod. More generally, it suffices that the product in the Koszul homology of $k[H]$ in homological degrees $1$ and $n-2$ vanishes. We also prove the conjecture when $H$ has embedding dimension four.
\end{abstract}
\maketitle
\section{Introduction}

Let $H=\langle a_1,\ldots,a_n\rangle$ be a numerical semigroup minimally generated by $a_1,\ldots,a_n$, where $n\ge3$, and let $k$ be a field. Thus $H$ has embedding dimension $n$. Set
\[
S=k[X_1,\ldots,X_n],\qquad R=k[H]\subseteq k[t],\qquad \deg X_i=a_i,
\]
and let
$\varphi:S\longrightarrow R$ be given by $\varphi(X_i)=t^{a_i}$ for $1\le i\le n$.
Put $I_H=\Ker\varphi$. We write $\n=(X_1,\ldots,X_n)$ and
$\widehat R=k[[H]]$. For a matrix $M$ over $S$, let $\detid_2(M)$ denote
the ideal generated by its $2\times2$ minors.

The pseudo-Frobenius numbers of $H$ are
\[
\PF(H)=\{f\in\mathbb Z\setminus H\mid f+a_i\in H
\text{ for every }1\le i\le n\}.
\]
The number $|\PF(H)|$ is called the \textit{type} of $H$, denoted by $\type(H)$.
Since the graded canonical module of $R$ is minimally generated by
$\{t^{-f} \mid f \in \PF(H)\} \subseteq k[t, t^{-1}]$
\cite{GotoWatanabe}, $\type(H)$ is equal to the Cohen--Macaulay type of
$R$.
In general, the set $\PF(H)$ does not determine $H$, even if the
embedding dimension is fixed. For example, the distinct numerical
semigroups
$H_1=\langle3,8,13\rangle$ and $H_2=\langle4,7,9\rangle$ both have
embedding dimension three and satisfy
$\PF(H_1)=\PF(H_2)=\{5,10\}$. Nevertheless, the conjecture below
predicts that, when $\PF(H)$ is an arithmetic progression of length
$n-1$, the ideal $I_H$ is generated by the $2\times2$ minors of a
$2\times n$ matrix whose entries are powers of the variables. In
particular, the total Betti numbers of $R$ over $S$ would then depend only on
$n$, through the Eagon--Northcott resolution.
The main object of this paper is the following conjecture, due to D.~T.~Cuong, D.~V.~Kien, N.~Matsuoka, and H.~L.~Truong.

\begin{conjecture}[Cuong--Kien--Matsuoka--Truong]\label{conj:main}
The following conditions are equivalent.
\begin{enumerate}
\item[\textup{(MP)}]\label{cond:homogeneous-matrix} The ideal $I_H$ has a \emph{maximal-minor presentation}: there are homogeneous elements $f_1,\ldots,f_n,g_1,\ldots,g_n\in S$ of positive degrees such that
\[
I_H=\detid_2\begin{pmatrix}
f_1&f_2&\cdots&f_n\\
g_1&g_2&\cdots&g_n
\end{pmatrix}.
\]
\item[\textup{(PP)}]\label{cond:pure-power-matrix} The ideal $I_H$ has a \emph{pure-power presentation}: after a permutation of $a_1,\ldots,a_n$, there are positive integers $\ell_1,\ldots,\ell_n,m_1,\ldots,m_n$ such that
\[
I_H=\detid_2\begin{pmatrix}
X_1^{\ell_1}&X_2^{\ell_2}&\cdots&X_n^{\ell_n}\\
X_2^{m_2}&X_3^{m_3}&\cdots&X_1^{m_1}
\end{pmatrix}.
\]
\item[\textup{(AP)}]\label{cond:AP} The set $\PF(H)$ is an \emph{arithmetic progression} of length $n-1$: there are $h,\alpha\in\mathbb Z$ with $\alpha>0$ such that
\[
\PF(H)=\{h+j\alpha\mid1\le j\le n-1\}.
\]
\end{enumerate}
\end{conjecture}

Conjecture~\ref{conj:main} arose from discussions among D.~T.~Cuong, D.~V.~Kien, N.~Matsuoka, and H.~L.~Truong and was first proposed in \cite[Remark~1]{KM}, with its present formulation given in \cite[Conjecture~1.1]{KMO}. The implication \textup{(PP)} $\Rightarrow$ \textup{(MP)} is immediate, and \textup{(MP)} $\Rightarrow$ \textup{(AP)} follows from the Eagon--Northcott resolution \cite[Theorem~3]{KM}. Thus the open part is \textup{(AP)} $\Rightarrow$ \textup{(PP)}. It is known for several classes, including maximal embedding dimension \cite{KM}, the case $h=0$ \cite{GKMT}, the case where $k[H]/(t^a)$ is stretched for
$a=\min(H\setminus\{0\})$ \cite[Theorem~1.2]{KMO}, and the range in which the multiplicity is at most $2n-2$ \cite{Takahashi}.

Kien proved the implication \textup{(AP)} $\Rightarrow$ \textup{(PP)} when $h$ has a unique factorization \cite[Theorem~3.7]{Kien}. Under condition~\textup{(AP)}, the unique factorization of $h$ is equivalent to that of $h+n\alpha$ \cite[Corollary~3.8]{Kien}. Murai and Takahashi \cite[Theorem~1.2]{MT} gave a numerical criterion for condition~\textup{(PP)} in arbitrary embedding dimension. They proved that the numerical identities expressing the pseudo-Frobenius numbers in terms of the exponents of the proposed matrix suffice to show that its maximal minors generate $I_H$.

Our first main result relates the three conditions to the Golod property.

\begin{theorem}\label{thm:main}
One has
\[
\textup{(MP)}
\quad\Longleftrightarrow\quad
\textup{(PP)}
\quad\Longleftrightarrow\quad
\textup{(AP)}\text{ and $\widehat R$ is Golod}.
\]
\end{theorem}

We obtain this theorem from the following stronger criterion involving multiplication in Koszul homology.

\begin{theorem}\label{thm:intro-tor}
Assume that, for some $h,\alpha\in\mathbb Z$ with $\alpha>0$,
\[
\PF(H)=\{h+j\alpha\mid1\le j\le n-1\},
\]
and
\[
\Tor_1^S(R,k)\,\Tor_{n-2}^S(R,k)=0.
\]
Then condition~\textup{(PP)} in Conjecture~\ref{conj:main} holds. In particular, under \textup{(AP)}, condition~\textup{(PP)} holds if $\widehat R$ is Golod.
\end{theorem}

Together with \textup{(PP)} $\Rightarrow$ \textup{(MP)} $\Rightarrow$ \textup{(AP)} and the fact that condition~\textup{(MP)} makes $\widehat R$ Golod, Theorem~\ref{thm:intro-tor} proves Theorem~\ref{thm:main} and reduces the open implication in Conjecture~\ref{conj:main} to a question about multiplication in Koszul homology.

In embedding dimension four, Murai and Takahashi \cite[Theorem~1.3]{MT} proved Conjecture~\ref{conj:main} under an additional assumption on an Ap\'ery set. More precisely, they assumed that, for some $1\le i\le4$, every element of
\[
\Ap(H,a_i)=\{w\in H\mid w-a_i\notin H\}
\]
has a unique expression as a non-negative integral combination of $a_1,a_2,a_3,a_4$. For example, $H=\langle12,14,15,17\rangle$ satisfies \textup{(AP)} but fails this Ap\'ery-set hypothesis, although Kien's criterion applies since $h=31$ has the unique factorization $14+17$. Our second main result removes this additional assumption.

\begin{theorem}\label{thm:intro-four}
If $H$ is minimally generated by four elements, then the three conditions in Conjecture~\ref{conj:main} are equivalent.
\end{theorem}

In the proofs of Theorems~\ref{thm:intro-tor} and~\ref{thm:intro-four}, we use Koszul homology to show that $h$ has a unique factorization and then apply Kien's theorem to obtain \textup{(PP)}. Under condition~\textup{(AP)}, the elements $h$ and $h+n\alpha$ belong to $H$ and extend the arithmetic progression by one term on either side.
The nonzero Koszul products used in the proof of
Theorem~\ref{thm:intro-tor} come from Poincar\'e duality for the
Artinian Gorenstein quotients of $R$ associated to $h$ and $h+n\alpha$. For Theorem~\ref{thm:intro-four}, the classification of Tor algebras of quotients of regular local rings of projective dimension three reduces the required relations to one binomial. The type bound for three-generated numerical semigroups then shows that the two terms of this binomial are powers of two variables.

In Section~\ref{sec:prelim}, we establish the general implications in Conjecture~\ref{conj:main} and recall the required facts on numerical semigroups and Koszul homology. In Section~\ref{sec:boundary}, we give a complete intersection criterion for the unique factorizations of $h$ and $h+n\alpha$. Section~\ref{sec:tor} is devoted to the proofs of Theorems~\ref{thm:main} and~\ref{thm:intro-tor}. In Section~\ref{sec:four}, we prove Theorem~\ref{thm:intro-four}.

\section{Preliminaries}\label{sec:prelim}

We first record the general implications in Conjecture~\ref{conj:main} and a consequence of condition~\textup{(MP)}. Put $Q=k[[X_1,\ldots,X_n]]$.

\begin{proposition}\label{prop:maximal-minor-golod}
Let $f_1,\ldots,f_n,g_1,\ldots,g_n$ belong to the maximal ideal of $Q$, and put
\[
M=
\begin{pmatrix}
f_1&f_2&\cdots&f_n\\
g_1&g_2&\cdots&g_n
\end{pmatrix}.
\]
If $\operatorname{ht}\detid_2(M)=n-1$, then $Q/\detid_2(M)$ is Golod.
\end{proposition}

\begin{proof}
Let
\[
T_{\mathrm{gr}}=k[Y_1,\ldots,Y_n,Z_1,\ldots,Z_n]/\detid_2
\begin{pmatrix}
Y_1&\cdots&Y_n\\
Z_1&\cdots&Z_n
\end{pmatrix}
\]
and
\[
T=k[[X_1,\ldots,X_n,Y_1,\ldots,Y_n,Z_1,\ldots,Z_n]]/\detid_2
\begin{pmatrix}
Y_1&\cdots&Y_n\\
Z_1&\cdots&Z_n
\end{pmatrix}.
\]
The ring $T_{\mathrm{gr}}$ is Golod by
\cite[Theorem~2.4\textup{(1)}]{BriggsChauDeStefani}; see also
\cite[Theorem~4]{HerzogReinerWelker}. The Golod property is
unchanged by completion; see the proof of
\cite[Proposition~5.2.4\textup{(4)}]{Avramov}. Adjoining the variables
$X_1,\ldots,X_n$ preserves the Golod property by
\cite[Proposition~5.2.4\textup{(2)}]{Avramov}, so $T$ is Golod. The
Eagon--Northcott resolution shows that $T$ is Cohen--Macaulay of dimension
$2n+1$.

The quotient of $T$ by
\[
Y_1-f_1,\ldots,Y_n-f_n,Z_1-g_1,\ldots,Z_n-g_n
\]
is isomorphic to $Q/\detid_2(M)$, which has dimension one by the height assumption. Hence
these $2n$ elements form a $T$-regular sequence. Their linear
parts are linearly independent, so each element remains outside the square
of the maximal ideal after the preceding elements have been killed.
Successive applications of \cite[Proposition~5.2.4\textup{(3)}]{Avramov}
now show that $Q/\detid_2(M)$ is Golod.
\end{proof}

Condition~\textup{(PP)} plainly implies condition~\textup{(MP)}. Under condition~\textup{(PP)}, each minor is a homogeneous binomial, so $\ell_i a_i-m_{i+1}a_{i+1}$ is independent of $i$, where $a_{n+1}=a_1$ and $m_{n+1}=m_1$. Thus \cite[Theorem~3]{KM} shows that condition~\textup{(AP)} holds. More generally, that theorem proves \textup{(MP)} $\Rightarrow$ \textup{(AP)} in Conjecture~\ref{conj:main}. Since $\operatorname{ht}I_H=n-1$, Proposition~\ref{prop:maximal-minor-golod} also shows that condition~\textup{(MP)} makes $\widehat R$ Golod.

We next recall some facts about numerical semigroups. We use the same notation as in the introduction.
Put $\sigma=\sum_{i=1}^n a_i$.
For $b,c\in\mathbb Z$, we write
\[
b\le_H c\quad\Longleftrightarrow\quad c-b\in H.
\]

We say that an element $s\in H$ has a \textit{unique factorization} if there is exactly one vector $(u_1,\ldots,u_n)\in\N^n$ with $s=\sum_{i=1}^n u_i a_i$.

\begin{lemma}\label{lem:semigroup}
Let $h,\alpha\in\mathbb Z$ with $\alpha>0$, and assume that
\[
\PF(H)=\{h+j\alpha\mid1\le j\le n-1\}.
\]
Put $h^+=h+n\alpha$.
Then the following assertions hold.
\begin{enumerate}
	\item $h,h^+ \in H$ and $j \alpha \notin H$ for every $1 \le j \le n-1$.
	\item For every $s \in H$,
	\[
	s+\alpha\notin H\quad\Longleftrightarrow\quad h-s\in H,
	\qquad
	s-\alpha\notin H\quad\Longleftrightarrow\quad h^+-s\in H.
	\]
\end{enumerate}
\end{lemma}
\begin{proof}
We first recall that, for every integer $z\notin H$, there exists $f\in\PF(H)$ such that $f-z\in H$. Indeed, we may take $f$ to be the largest integer in $(z+H)\setminus H$.

(1) Since $h+\alpha,h+2\alpha\in\PF(H)$, we have $\alpha\notin H$.
Applying the property recalled above to $\alpha\notin H$, we have
$h+(p-1)\alpha\in H$ for some $1\le p\le n-1$.
If $p\ge2$, then $h+(p-1)\alpha\in\PF(H)$, which is impossible. Hence $p=1$ and $h\in H$; see also \cite[Lemma~2]{KM}.
Now, for $1 \le j \le n-1$, $j\alpha\in H$ would imply $h+j\alpha\in H$, a contradiction.
Hence we have $j\alpha \notin H$ for every $1 \le j \le n-1$.
Moreover, $h+(n-1)\alpha$ is the Frobenius number of $H$, so $h^+\in H$.

(2) The first equivalence is proved in \cite[Lemma~3.3]{Kien}. The second follows by the same argument, which we include for completeness. Let $s\in H$ and suppose that $s-\alpha\notin H$.
By the property recalled at the beginning of the proof, $h+(p+1)\alpha-s\in H$ for some $1\le p\le n-1$.
Adding $s$, we obtain $h+(p+1)\alpha\in H$, which forces $p=n-1$. Hence $h^+-s\in H$.
Conversely, since $(s-\alpha)+(h^+-s)=h+(n-1)\alpha\notin H$, we have $s-\alpha\notin H$ whenever $h^+-s\in H$.
\end{proof}
We now fix notation for the Koszul complex and its multiplication.

Let $F=\bigoplus_{j=1}^nR(-a_j)$ be a graded free $R$-module with basis $e_1,\ldots,e_n$, where $\deg e_j=a_j$ for $1 \le j \le n$. For a graded $R$-module $M$, we put $K_\bullet(M)=K_\bullet(t^{a_1},\ldots,t^{a_n};M)$, the Koszul complex of $M$ with respect to $t^{a_1},\ldots,t^{a_n}$. Thus
\[
K_i(M)=\bigwedge_R^i F\otimes_R M,\qquad \hml_i(K_\bullet(M))=\Tor_i^S(M,k).
\]
For $\Lambda=\{j_1<\cdots<j_i\}\subseteq\{1,\ldots,n\}$, we put $e_{\Lambda}=e_{j_1}\wedge\cdots\wedge e_{j_i}$ and $a_{\Lambda}=\sum_{j\in \Lambda}a_j$. The differential is given by
\[
\partial_i(me_{\Lambda})=\sum_{r=1}^i(-1)^{r-1}t^{a_{j_r}}m e_{\Lambda\setminus\{j_r\}}
\]
for $m \in M$.
We define the multiplication $K_i(R) \times K_j(R) \to K_{i+j}(R)$ by 
\[
(t^s e_{\Lambda})(t^{s'}e_{\Gamma})=t^{s+s'}(e_{\Lambda}\wedge e_{\Gamma}),
\]
where $s,s'\in H$, $i,j \ge 0$, and $\Lambda,\Gamma\subseteq\{1,\ldots,n\}$ with $|\Lambda| = i$ and $|\Gamma| = j$. 
Extending this formula $k$-bilinearly gives multiplication in $K_\bullet(R)$.

Finally, we recall the following result of Kien. It gives condition~\textup{(PP)} once $h$ is known to have a unique factorization and will be applied in the proof of Theorem~\ref{thm:tor}.

\begin{theorem}[{\cite[Theorem~3.7]{Kien}}]\label{thm:kien}
Let $h,\alpha\in\mathbb Z$ with $\alpha>0$, and assume that
\[
\PF(H)=\{h+j\alpha\mid1\le j\le n-1\}.
\]
Suppose that $h$ has the unique factorization
$h=\sum_{i=1}^n u_i a_i$. Then, after a permutation, there are positive
integers $v_1,\ldots,v_n$ such that
\[
(u_i+1)a_i+\alpha=v_{i+1}a_{i+1}\quad(1\le i\le n),\qquad
I_H=\detid_2\begin{pmatrix}X_1^{u_1+1}&\cdots &X_{n-1}^{u_{n-1}+1} & X_n^{u_n+1}\\X_2^{v_2}&\cdots & X_n^{v_n} &X_1^{v_1}\end{pmatrix}.
\]
Here the indices are read modulo $n$.
\end{theorem}

\begin{corollary}[{\cite[Corollary~3.8]{Kien}}]\label{cor:endpoints-uf}
Assume \textup{(AP)} and put $h^+=h+n\alpha$. Then $h$ has a unique
factorization if and only if $h^+$ has a unique factorization.
 If these equivalent conditions hold, then
condition~\textup{(PP)} in Conjecture~\ref{conj:main} holds.
\end{corollary}

\section{Complete intersections and unique factorizations}\label{sec:boundary}

Throughout this section, $H=\langle a_1,\ldots,a_n\rangle$, $S$, $R$,
$\varphi:S\to R$, $I_H$, $\n$, $\sigma$, and $\le_H$ retain the meanings
fixed above. Thus $\varphi(X_i)=t^{a_i}$ for every $1\le i\le n$. Let
$h,\alpha\in\mathbb Z$ with $\alpha>0$, and assume that
\[
\PF(H)=\{h+j\alpha\mid1\le j\le n-1\}.
\]
We retain the notation $h^+=h+n\alpha$. We call $h$ and $h^+$ the endpoints of the extended arithmetic progression
$h,\ h+\alpha,\ldots,h+(n-1)\alpha,\ h+n\alpha$.
\begin{lemma}\label{lem:boundary-identities}
The following assertions hold true.
\begin{enumerate}
	\item Let $s\in H$ and $1\le q\le n-1$. If $s+q\alpha\in H$, then $s+j\alpha\in H$ for every $0\le j\le q$.
	\item Let $d\in H$ with $d\le_Hh$ and $d\le_Hh^+$. Then, for every $b\in H$ and $1\le r\le n-1$,
\[
d-b\in H\iff h-d+b+r\alpha\notin H.
\]
Furthermore, let $1\le i<j\le n$. If $d-a_i,d-a_j\in H$ but
$d-a_i-a_j\notin H$, then
\[
h-d+a_i+a_j+r\alpha\in\langle a_\ell\mid\ell\ne i,j\rangle
\]
for every $1\le r\le n-1$.
\end{enumerate}
\end{lemma}

\begin{proof}
\textup{(1)} Suppose that there exists $0 \le j \le q$ such that $s + j\alpha \notin H$, and let $j$ be the least such integer.
Then $1\le j<q$, and the first equivalence in Lemma~\ref{lem:semigroup}~\textup{(2)} gives $h-s-(j-1)\alpha\in H$. Adding $s+q\alpha$ yields $h+(q-j+1)\alpha\in H$, contrary to the pseudo-Frobenius hypothesis.

\textup{(2)} Let $b\in H$ and $1\le r\le n-1$.
If $d-b\in H$ and $h-d+b+r\alpha\in H$, then $h+r\alpha = (d-b) + (h-d+b+r\alpha) \in H$, a contradiction. Conversely, suppose that $d-b\notin H$.
The first equivalence in Lemma~\ref{lem:semigroup}~\textup{(2)} gives $h-d+b+\alpha\in H$. Since
$h-d+b+n\alpha=(h^+-d)+b\in H$, the second equivalence in Lemma~\ref{lem:semigroup}~\textup{(2)} gives
$h-d+b+(n-1)\alpha\in H$. Applying \textup{(1)} with
$s=h-d+b+\alpha$ and $q=n-2$ gives
$h-d+b+r\alpha\in H$ for every $1\le r\le n-1$.

For the additional assertion, let $1\le r\le n-1$. Applying the equivalence
in \textup{(2)} to $b=a_i+a_j$ gives
$
h-d+a_i+a_j+r\alpha\in H.
$
Applying the equivalence in \textup{(2)} to $b=a_i$ gives
$
h-d+a_i+r\alpha\notin H.
$
The case $b=a_j$ gives $h-d+a_j+r\alpha\notin H$. If $a_i$ occurs in a factorization of
$h-d+a_i+a_j+r\alpha$, then $h-d+a_j+r\alpha\in H$, a contradiction.
Interchanging $i$ and $j$ shows that $a_j$ does not occur in any such
factorization. Therefore
\[
h-d+a_i+a_j+r\alpha\in\langle a_\ell\mid\ell\ne i,j\rangle.
\]
\end{proof}

For $x\in H$, set
\[
D_x=\{b\in H\mid x-b\in H\},\qquad
\mathfrak a_x=(t^b\mid b\in H\setminus D_x),
\]
and define $B_x=R/\mathfrak a_x$. Let
\[
\rho_x:S\longrightarrow B_x,
\qquad
\rho_x(X_i)=\overline{t^{a_i}},
\]
and put $J_x=\Ker\rho_x$. Since $\rho_x$ is the composite of
$\varphi:S\to R$ and the quotient map $R\to R/\mathfrak a_x$, we have
\[
J_x=\varphi^{-1}(\mathfrak a_x).
\]
The set $D_x$ is finite and downward closed with respect to $\le_H$. 
The algebra $B_x$ has $k$-basis $\{\overline{t^b}\mid b\in D_x\}$ and socle $k\overline{t^x}$. Hence it is Artinian Gorenstein. 
Lemma~\ref{lem:semigroup}~\textup{(2)} gives
\[
\mathfrak a_h=R:_R t^\alpha,
\qquad
\mathfrak a_{h^+}=t^\alpha(R:_R t^\alpha).
\]
Furthermore, let
\[
V_x=\bigoplus_{d\le_Hx}[I_H/\n I_H]_d.
\]

\begin{lemma}\label{lem:divisor-algebras}
For $x \in H$, the following conditions are equivalent: 
\begin{enumerate}
	\item $x$ has a unique factorization.
	\item $J_x$ is a monomial ideal.
	\item $V_x=0$.
\end{enumerate} 
If $x=\sum_{i=1}^n u_i a_i$ is the unique factorization of $x$, then
\[
J_x=(X_1^{u_1+1},\ldots,X_n^{u_n+1}).
\]
Furthermore, suppose that the class of a homogeneous binomial $f$ belonging
to a minimal homogeneous generating set of $I_H$ spans $V_x$ as a
$k$-vector space; that is, $V_x=k\overline f$. Then, for each $y\in H$ with
$y\le_Hx$,
$[I_H]_y=[(f)]_y$.
In particular, there are exactly two monomials in $S$ of degree $\deg f$,
and they are relatively prime.
\end{lemma}

\begin{proof}
Suppose that $u=(u_1,\ldots,u_n)$ is the unique factorization of $x$. The image of $X^v$ in $B_x$ is nonzero precisely when $v$ extends to a factorization of $x$, that is, when $v\le u$ coordinatewise.
If $v,w\le u$ and $X^v$ and $X^w$ have the same degree, then $u-v+w$ is a factorization of $x$, so $v=w$. Thus the monomials with nonzero image have distinct degrees, and $J_x=(X_1^{u_1+1},\ldots,X_n^{u_n+1})$.
Conversely, if $u$ and $v$ are distinct factorizations of $x$, then $X^u-X^v\in J_x$ while neither term belongs to $J_x$. This proves \textup{(1)}~$\Leftrightarrow$~\textup{(2)}.

For $y\le_Hx$, only minimal homogeneous generators of $I_H$ of degrees $d\le_Hy$ contribute to $[I_H]_y$. Conversely, each generator of degree $d\le_Hx$ has a nonzero monomial multiple of degree $x$. Hence $V_x=0$ if and only if $[I_H]_x=0$, which is equivalent to \textup{(1)}.
If $V_x=k\overline f$, choose a minimal homogeneous generating set containing $f$. Then $f$ is its only member of degree $\le_Hx$, so $[I_H]_y=[(f)]_y$ for every $y\le_Hx$.
In particular, $[I_H]_{\deg f}=kf$. Since all monomials of degree $\deg f$ map to $t^{\deg f}$, there are exactly two. Their relative primality follows from the minimality of $f$.
\end{proof}

\begin{remark}
By \cite[Lemma~2.4(3)]{EtoWatanabe}, $J_x$ is the annihilator in $S$ of the inverse polynomial associated to $x$.
\end{remark}

We now prove the complete intersection criterion used in the four-generated case.

\begin{theorem}\label{thm:boundary-reconstruction}
If both $J_h$ and $J_{h^+}$ are complete intersections, then $h$ and $h^+$
have unique factorizations.
\end{theorem}
\begin{proof}
Let $L=R+Rt^\alpha\subseteq k[t,t^{-1}]$, regarded as a graded $S$-module
through the homomorphism $\varphi:S\to R$. Let $Se_0\oplus Se_1$ be the
graded free $S$-module with basis $e_0,e_1$, where $\deg e_0=0$ and
$\deg e_1=\alpha$. Define the graded $S$-linear map
\[
\psi:Se_0\oplus Se_1\longrightarrow L,
\qquad
\psi(pe_0+qe_1)=\varphi(p)+\varphi(q)t^\alpha
\]
for $p,q\in S$. Put $N=\Ker\psi$. Since $\varphi$ is
surjective, so is $\psi$, and
we have a graded presentation
\[
0\longrightarrow N\longrightarrow Se_0\oplus Se_1
\xrightarrow{\ \psi\ }L\longrightarrow0.
\]
This presentation is minimal because $\alpha\notin H$. For $p,q\in S$,
we write $(p,q)=pe_0+qe_1$ and define
\[
\pi_0:N\longrightarrow Se_0,\quad \pi_0(p,q)=pe_0,
\qquad
\pi_1:N\longrightarrow Se_1,\quad \pi_1(p,q)=qe_1.
\]
Then the equalities
\[
\mathfrak a_h=R:_R t^\alpha,
\qquad
\mathfrak a_{h^+}=t^\alpha(R:_R t^\alpha)
\]
give exact sequences
\[
0\longrightarrow I_He_0\longrightarrow N
\xrightarrow{\ \pi_1\ }J_he_1\longrightarrow0,
\qquad
0\longrightarrow I_He_1\longrightarrow N
\xrightarrow{\ \pi_0\ }J_{h^+}e_0\longrightarrow0.
\]
We verify the first exact sequence; the second is analogous. Let
$(p,q)\in N$. Then
\[
\varphi(q)t^\alpha=-\varphi(p)\in R.
\]
Hence $\varphi(q)\in R:_Rt^\alpha=\mathfrak a_h$, and therefore
$q\in J_h$. Conversely, if $q\in J_h$, then
$\varphi(q)t^\alpha\in R$. Since $\varphi$ is surjective, we may choose $p\in S$ such that
$\varphi(p)=-\varphi(q)t^\alpha$, and then $(p,q)\in N$. Thus the image
of $\pi_1$ is $J_he_1$. Its kernel consists of
the elements $(p,0)$ with $\varphi(p)=0$, and hence is $I_He_0$.
The same argument applies to $\pi_0$. The equality
$\mathfrak a_{h^+}=t^\alpha(R:_Rt^\alpha)$ shows that the image of
$\pi_0$ is $J_{h^+}e_0$, and its kernel is $I_He_1$.

The inclusion $N\to Se_0\oplus Se_1$ induces a map
\[
\bigwedge^2N\longrightarrow\bigwedge^2(Se_0\oplus Se_1)\cong S.
\]
Let $\fkc$ be the image of this map. Suppose that
$N=Sv_1+\cdots+Sv_\ell$, where $v_i=(r_i,s_i)$ and $r_i,s_i\in S$
for $1\le i\le\ell$. Since $\bigwedge^2N$ is generated by
$v_i\wedge v_j$, where $1\le i<j\le\ell$, and the image of
$v_i\wedge v_j$ is $r_is_j-r_js_i$, we have
$\fkc=(r_is_j-r_js_i\mid 1\le i<j\le\ell)$.
Moreover, for $1\le i<j\le\ell$, the equalities
\[
\varphi(r_i)=-\varphi(s_i)t^\alpha,
\qquad
\varphi(r_j)=-\varphi(s_j)t^\alpha,
\]
give $\varphi(r_is_j-r_js_i)=0$. Hence $\fkc\subseteq I_H$.

Set
$U=I_H/(\n I_H+\fkc)$, $\overline N=N/\n N$
and define two $S$-linear maps
\[
\iota_0:I_H\longrightarrow\overline N,\quad
\iota_0(f)=\overline{(f,0)},
\qquad
\iota_1:I_H\longrightarrow\overline N,\quad
\iota_1(f)=\overline{(0,f)}.
\]
We claim that both maps have kernel $\n I_H+\fkc$. If $(p,q)\in N$,
then $p\in J_{h^+}$ and $q\in J_h$. In particular, $p,q\in\n$.
For $1\le i<j\le\ell$, we have
\[
s_jv_i-s_iv_j
=s_j(r_i,s_i)-s_i(r_j,s_j)
=(r_is_j-r_js_i,0)
\]
and
\[
r_iv_j-r_jv_i
=r_i(r_j,s_j)-r_j(r_i,s_i)
=(0,r_is_j-r_js_i).
\]
Since $r_i,r_j,s_i,s_j\in\n$, the left-hand sides belong to $\n N$.
Thus $(r_is_j-r_js_i,0)$ and $(0,r_is_j-r_js_i)$ belong to $\n N$.
Since the elements
$r_is_j-r_js_i$, where $1\le i<j\le\ell$, generate $\fkc$, it
follows that $\fkc$ is contained in the kernel of each map. Moreover,
if $f\in I_H$, then $(f,0),(0,f)\in N$. Thus $\n I_H$ is also
contained in the kernel of each map. Therefore
$\n I_H+\fkc$ is contained in both kernels.

Conversely, consider $\iota_0$. Choose a homogeneous regular sequence
$g_1,\ldots,g_n$ minimally generating $J_h$. For $1\le i\le n$,
choose $f_i\in S$ such that $(f_i,g_i)\in N$. Since
$\Ker\pi_1=I_He_0$, if $\gamma_1,\ldots,\gamma_s$ generate $I_H$, then
\[
(f_1,g_1),\ldots,(f_n,g_n),(\gamma_1,0),\ldots,(\gamma_s,0)
\]
generate $N$. If $(f,0)\in\n N$, then there are
$\lambda_1,\ldots,\lambda_n,\mu_1,\ldots,\mu_s\in\n$ such that
\[
(f,0)=\sum_{i=1}^n \lambda_i(f_i,g_i)
+\sum_{j=1}^s \mu_j(\gamma_j,0).
\]
Applying $\pi_1$ to this equality, we obtain
$\sum_{i=1}^n \lambda_ig_i=0$. Since $g_1,\ldots,g_n$ is a regular
sequence, its first syzygies are generated by the Koszul relations.
Thus there are $c_{ij}\in S$, where $1\le i<j\le n$, such that
\[
\lambda_i=\sum_{j=i+1}^n c_{ij}g_j-\sum_{j=1}^{i-1}c_{ji}g_j
\qquad(1\le i\le n).
\]
Therefore
\[
\sum_{i=1}^n \lambda_if_i
=\sum_{1\le i<j\le n}c_{ij}(g_jf_i-g_if_j)\in \fkc.
\]
Applying $\pi_0$ to the equality defining $(f,0)$ gives
$f=\sum_{i=1}^n \lambda_if_i+\sum_{j=1}^s \mu_j\gamma_j$.
The first sum belongs to $\fkc$. The second belongs to $\n I_H$
because $\mu_j\in\n$ and $\gamma_j\in I_H$ for $1\le j\le s$. Hence
$f\in\fkc+\n I_H$. The same argument applied to $\iota_1$, using a
regular sequence minimally generating $J_{h^+}$, proves that
$\Ker\iota_1=\n I_H+\fkc$.

For a graded vector space $G$, write $G(\alpha)$ for the shift defined by
$[G(\alpha)]_d=[G]_{d+\alpha}$. Passing the two exact sequences with middle
term $N$ modulo $\n$ and using the kernel computations gives
\[
0\longrightarrow U\longrightarrow \overline N\longrightarrow (J_h/\n J_h)(-\alpha)\longrightarrow0,
\qquad
0\longrightarrow U(-\alpha)\longrightarrow \overline N\longrightarrow J_{h^+}/\n J_{h^+}\longrightarrow0.
\]
For a finite-dimensional graded vector space $G$, put $\delta(G)=\sum_{d\in\mathbb Z} d\dim_k[G]_d$. Since $J_h$ and $J_{h^+}$ are $\n$-primary complete intersections, $\dim_k(J_h/\n J_h)=\dim_k(J_{h^+}/\n J_{h^+})=n$. Taking degree sums gives
\[
\delta(\overline N)=\delta(U)+\delta(J_h/\n J_h)+n\alpha, \qquad \delta(\overline N)=\delta(U)+(\dim_kU)\alpha
+\delta(J_{h^+}/\n J_{h^+})
\]
from the exact sequences above. Therefore
\[
\delta(J_{h^+}/\n J_{h^+})-\delta(J_h/\n J_h)
=(n-\dim_kU)\alpha.
\]
For an Artinian weighted complete intersection, the socle degree is the sum
of the relation degrees minus $\sigma = \sum_{i=1}^n a_i$; see
\cite[Corollary~3.6.14]{BrunsHerzog}.
Since the socle degrees of $B_h$ and $B_{h^+}$ are $h$ and $h^+$, respectively, we have
\[
\delta(J_h/\n J_h)-\sigma=h,
\qquad
\delta(J_{h^+}/\n J_{h^+})-\sigma=h^+.
\]
Thus
\[
\delta(J_{h^+}/\n J_{h^+})-\delta(J_h/\n J_h)=h^+-h=n\alpha,
\]
and consequently $\dim_k U = 0$, whence $U=0$. The two exact sequences with middle term $\overline N$ now give
\[
\overline N\cong (J_h/\n J_h)(-\alpha)
\cong J_{h^+}/\n J_{h^+},
\]
and hence $\dim_k\overline N=n$.

For every $d\in\mathbb Z$, the set
\[
\{me_0\mid m\text{ is a monomial of degree }d\}
\cup
\{me_1\mid m\text{ is a monomial of degree }d-\alpha\}
\]
is a $k$-basis of $[Se_0\oplus Se_1]_d$. Under $\psi$, each element of
this basis maps to $t^d$. Hence
$[N]_d$ is spanned by elements of the following three forms:
\[
(m-m',0),\qquad (0,m-m'),\qquad (m,-m'),
\]
where $\deg m=\deg m'=d$ in the first form,
$\deg m=\deg m'=d-\alpha$ in the second form, and
$\deg m=d$ and $\deg m'=d-\alpha$ in the third form. An element of either of the first two forms
has zero image in $\overline N$. In fact, $m-m'\in I_H$, and $U=0$
implies that both $\iota_0(m-m')$ and $\iota_1(m-m')$ are zero.
Since $\dim_k\overline N=n$, we may therefore choose the $n$ minimal
generators of $N$ as
$v_i=(p_i,-q_i)$ for $1\le i\le n$,
where $p_i$ and $q_i$ are monomials and
$\deg p_i-\deg q_i=\alpha$ for $1\le i\le n$.
Thus $p_1,\ldots,p_n$ and $q_1,\ldots,q_n$ are minimal generating
sets of $J_{h^+}$ and $J_h$, respectively. Each of these ideals is an
$\n$-primary monomial ideal with exactly $n$ generators, so its minimal
monomial generating set consists of one pure power of each variable.
Lemma~\ref{lem:divisor-algebras} now shows that $h$ and $h^+$ have unique
factorizations.
\end{proof}

\begin{remark}
The converse of Theorem~\ref{thm:boundary-reconstruction} follows from Lemma~\ref{lem:divisor-algebras}. For each $x\in\{h,h^+\}$, if $x$ has a unique factorization, then $J_x$ is generated by one pure power of each variable and hence is a complete intersection.
\end{remark}

\section{Koszul homology and the Golod case}\label{sec:tor}

For $0\le i\le n$, put
\[
A_i=\hml_i(K_\bullet(R))=\Tor_i^S(R,k),
\qquad
A=\bigoplus_{i=0}^nA_i.
\]
The multiplication defined on $K_\bullet(R)$ in Section~\ref{sec:prelim} makes $A$ a $k$-algebra. We identify $A_1$ with $I_H/\n I_H$ and write $[A_i]_d$ for the degree-$d$ component of $A_i$.

We first prove that the vanishing of $A_1A_{n-2}$ implies that $h$ has a unique factorization. We then apply Theorem~\ref{thm:kien} and deduce the Golod case of Conjecture~\ref{conj:main}.

\begin{proposition}\label{prop:boundary-products}
Let $d\in\mathbb Z$. Assume that, for some $h,\alpha\in\mathbb Z$ with $\alpha>0$,
\[
\PF(H)=\{h+j\alpha\mid1\le j\le n-1\}.
\]
If $d\le_Hh$, then, for every nonzero
$\xi\in[A_1]_d$,
there exists
\[
\eta\in[A_{n-2}]_{\sigma+h+\alpha-d}
\quad\text{such that}\quad \xi\eta\ne0.
\]
If $d\le_Hh^+$, then, for every nonzero
$\xi\in[A_1]_d$,
there exists
\[
\eta\in[A_{n-2}]_{\sigma+h+(n-1)\alpha-d}
\quad\text{such that}\quad \xi\eta\ne0.
\]
\end{proposition}

\begin{proof}
For the first assertion put $(x,\varepsilon)=(h,\alpha)$, and for the
second put $(x,\varepsilon)=(h^+,-\alpha)$. Write $f=x+\varepsilon$.
Then $f\in\PF(H)$, and Lemma~\ref{lem:semigroup}~\textup{(2)}
gives $\mathfrak a_x=R:_R t^\varepsilon$. Suppose that $d\le_H x$.
Since $d-a_\Lambda\le_H d\le_H x$ for every $\Lambda \subseteq \{1,\dots,n\}$, we have
$[\mathfrak a_x]_{d-a_\Lambda}=0$ and hence $[K_\bullet(\mathfrak a_x)]_d=0$.
The short exact sequence $0\to\mathfrak a_x\to R\to B_x\to0$
therefore induces an isomorphism $[K_\bullet(R)]_d\cong[K_\bullet(B_x)]_d$.

Put $e_{\{1,\ldots,n\}}=e_1\wedge\cdots\wedge e_n$, so that $\deg e_{\{1,\ldots,n\}}=\sigma=\sum_{j=1}^n a_j$.
For every $b\in\mathbb Z$, multiplication in $B_x$ induces a perfect pairing
\[
  [B_x]_b \times [B_x]_{x-b} \longrightarrow [B_x]_x.
\]
Indeed, if $b,x-b\in H$, the two factors are spanned by
$\overline{t^b}$ and $\overline{t^{x-b}}$, whose product spans $[B_x]_x$.
Otherwise both factors vanish.
Together with the exterior product, the pairing on $B_x$ gives
a perfect pairing
\[
  [K_i(B_x)]_b \times [K_{n-i}(B_x)]_{\sigma+x-b}
  \longrightarrow [K_n(B_x)]_{\sigma+x}
  = k\,\overline{t^x}e_{\{1,\ldots,n\}}
\]
for every $0\le i\le n$ and $b\in\mathbb Z$.
Since $\overline{t^x}$ lies in the socle of $B_x$ and $K_{n+1}(B_x)=0$,
the class of $\overline{t^x}e_{\{1,\ldots,n\}}$ spans the one-dimensional
space $[\hml_n(K_\bullet(B_x))]_{\sigma+x}$.
Since $K_{n+1}(B_x)=0$, the Leibniz rule shows that the differentials
are adjoint up to sign under these pairings. Thus the orthogonal complement
of the boundaries is the space of cycles in the complementary degree,
and the induced pairing on homology is perfect.

Let $z\in[K_1(R)]_d$ represent $\xi\ne0$. Its image $\bar z$ in
$[K_1(B_x)]_d$ represents a nonzero class by the isomorphism above.
By the perfect pairing on homology, we may choose a cycle
$v\in[K_{n-1}(B_x)]_{\sigma+x-d}$ such that
$[\bar z][\,v]=[\overline{t^x}e_{\{1,\ldots,n\}}]$.
Since $K_{n+1}(B_x)=0$, this gives
$\bar z\,v=\overline{t^x}e_{\{1,\ldots,n\}}$. Take a homogeneous lift $\widetilde v$ to
$K_{n-1}(R)$. 
Since $v$ is a cycle, the coefficients of $\partial\widetilde v$
belong to $\mathfrak a_x$. As $t^\varepsilon\mathfrak a_x\subseteq R$,
we obtain
\[
  w=t^\varepsilon\partial\widetilde v
    \in[K_{n-2}(R)]_{\sigma+f-d}.
\]
Moreover, $\partial w=t^\varepsilon\partial^2\widetilde v=0$.
Since $[R]_x\to[B_x]_x$ is injective,
$z\widetilde v=t^xe_{\{1,\ldots,n\}}$. 
Computing in $K_\bullet(k[t,t^{-1}])$ and using $\partial z=0$, we get
\[
 zw=-t^\varepsilon\partial(z\widetilde v)
    =-\sum_{i=1}^n(-1)^{i-1}t^{f+a_i}
         e_{\{1,\ldots,n\}\setminus\{i\}} .
\]
Since $f\in\PF(H)$, each coefficient belongs to $R$.
The terms involve distinct basis elements, so $zw\ne0$.
It represents a nonzero class since
$[K_n(R)]_{\sigma+f}=[R]_f e_{\{1,\ldots,n\}}=0$. Thus $\eta=[w]$ has the required
degree and $\xi\eta\ne0$.
\end{proof}

\Needspace{6\baselineskip}
We now obtain the following theorem.

\begin{theorem}\label{thm:tor}
Suppose that, for some $h,\alpha\in\mathbb Z$ with $\alpha>0$,
\[
\PF(H)=\{h+j\alpha\mid1\le j\le n-1\}.
\]
Assume also that $A_1A_{n-2}=0$. Then condition~\textup{(PP)} in Conjecture~\ref{conj:main} holds.
\end{theorem}

\begin{proof}
Suppose that $h$ does not have a unique factorization. By Lemma~\ref{lem:divisor-algebras}, we have $V_h\ne0$. Thus we may choose $d\in\mathbb Z$ and $0\ne\xi\in[A_1]_d$ such that $h-d\in H$.
Proposition~\ref{prop:boundary-products} gives
\[
\eta\in[A_{n-2}]_{\sigma+h+\alpha-d}
\quad\text{such that}\quad \xi\eta\ne0,
\]
contrary to $A_1A_{n-2}=0$. Hence $h$ has a unique factorization, and Theorem~\ref{thm:kien} gives condition~\textup{(PP)}.
\end{proof}

The Golod case follows as a corollary.

\begin{corollary}\label{cor:golod}
Assume that $\widehat R$ is Golod and that, for some $h,\alpha\in\mathbb Z$ with $\alpha>0$,
\[
\PF(H)=\{h+j\alpha\mid1\le j\le n-1\}.
\]
Then condition~\textup{(PP)} in Conjecture~\ref{conj:main} holds.
\end{corollary}
\begin{proof}
Since $\widehat R$ is Golod, all products of positive Koszul homology classes are zero by \cite[Remarks~5.2.1 and~5.2.3]{Avramov}.
The Koszul homology modules of $R$ are annihilated by its homogeneous maximal ideal. Therefore localization at this ideal and completion preserve these modules and their products.
Since $n\ge3$, we have $A_1A_{n-2}=0$.
Theorem~\ref{thm:tor} now gives condition~\textup{(PP)}.
\end{proof}

Theorem~\ref{thm:tor} and Corollary~\ref{cor:golod} together prove Theorem~\ref{thm:intro-tor}.
Condition~\textup{(PP)} implies condition~\textup{(MP)}. By Section~\ref{sec:prelim}, condition~\textup{(MP)} implies condition~\textup{(AP)} and makes $\widehat R$ Golod. Conversely, Corollary~\ref{cor:golod} gives condition~\textup{(PP)} when condition~\textup{(AP)} holds and $\widehat R$ is Golod. This proves Theorem~\ref{thm:main}.
In particular, we have the following.

\begin{corollary}\label{cor:two-presentations}
Conditions~\textup{(MP)} and~\textup{(PP)} in Conjecture~\ref{conj:main} are equivalent.
\end{corollary}

\begin{remark}
Under the hypotheses of \cite[Theorem~1.2]{KMO}, condition~(3) of that theorem implies \textup{(PP)} by Corollary~\ref{cor:golod}. To verify this, put $a=\min(H\setminus\{0\})$ and $B=\widehat R/(t^a)$, and assume that the square of the maximal ideal of $B$ is principal and that $\PF(H)$ has the form in Corollary~\ref{cor:golod}. Thus $B$ is stretched. The Cohen--Macaulay type and the embedding dimension of $B$ are both $n-1$. Hence $B$ is Golod by \cite[Theorem~5.4\textup{(2)(a)}]{CrollEtAl}. Since $t^a$ is a regular element of $\widehat R$ outside the square of its maximal ideal, $\widehat R$ is also Golod by \cite[Proposition~5.2.4\textup{(2)}]{Avramov}. We can therefore apply Corollary~\ref{cor:golod}.
\end{remark}

Under condition~\textup{(AP)}, Theorems~\ref{thm:main} and~\ref{thm:tor}, together with the vanishing of products of positive Koszul homology classes for Golod rings, show that each of conditions~\textup{(MP)} and~\textup{(PP)} is equivalent both to the Golod property of $\widehat R$ and to the equality $A_1A_{n-2}=0$.
Thus the open implication \textup{(AP)} $\Rightarrow$ \textup{(PP)} is equivalent to an affirmative answer to the following question.

\begin{question}\label{que:tor-vanishing}
Assume that, for some $h,\alpha\in\mathbb Z$ with $\alpha>0$,
\[
\PF(H)=\{h+j\alpha\mid1\le j\le n-1\}.
\]
Is it true that $A_1A_{n-2}=0$?
\end{question}

\section{The four-generated case}\label{sec:four}

Throughout this section, assume that $n=4$ and that there are $h,\alpha\in\mathbb Z$ with $\alpha>0$ such that
\[
\PF(H)=\{h+\alpha,h+2\alpha,h+3\alpha\}.
\]
The symbols $S$, $R$, $I_H$, $\n$, $h^+$, $A_i$, and $A$ retain the
meanings fixed in the preceding sections.

The Tor algebras of quotients $B=Q/J$ with $Q$ regular local and $\operatorname{pd}_QB=3$ are divided into several classes by Avramov, Kustin, and Miller \cite{AKM}; see also \cite{ChristensenVeliche}. We first record the consequence of this classification that we need.

\begin{proposition}\label{prop:pd-three}
Let $(Q,\mathfrak q,k)$ be a regular local ring and $B=Q/J$, where $J\subseteq\mathfrak q^2$ and $\operatorname{pd}_QB=3$. For $0\le i\le3$, put
\[
T_i=\Tor_i^Q(B,k),
\qquad
T=\bigoplus_{i=0}^3T_i.
\]
If $q=\dim_k(T_1T_2)\ge2$, then
$\dim_k\bigl(T_1/\Ann_{T_1}(T_2)\bigr)=1$.
\end{proposition}

\begin{proof}
Put $p=\dim_k(T_1T_1)$. If $B$ were Gorenstein, then $\dim_k T_3=1$, contradicting $q\ge2$ since $T_1T_2\subseteq T_3$. Since $q\ge2$, the classification places $T$ in class $\mathbf H(p,q)$ \cite[Proposition~2]{ChristensenVeliche}. For an algebra in this class, there are bases $\{e_i\}_{1\le i\le\dim_kT_1}$ of $T_1$, $\{b_j\}_{1\le j\le\dim_kT_2}$ of $T_2$, and $\{c_r\}_{1\le r\le\dim_kT_3}$ of $T_3$ such that the only nonzero products among the $e_i b_j$ are
\[
e_{p+1}b_{p+r}=c_r\qquad(1\le r\le q).
\]
Thus
\[
e_i b_j=0
\quad\text{if }(i,j)\ne(p+1,p+r)\text{ for every }1\le r\le q
\]
\cite[3.9]{BriggsEtAl}. Thus an element of $T_1$ annihilates $T_2$ exactly when its coefficient of $e_{p+1}$ is zero.
\end{proof}

\begin{remark}
Let
$\mu:T_1\longrightarrow\Hom_k(T_2,T_3), \qquad \mu(x)(y)=xy$.
Then $\Ker\mu=\Ann_{T_1}(T_2)$, so Proposition~\ref{prop:pd-three} says that $\dim_k\im\mu=1$. On the other hand, $q=\dim_k(T_1T_2)$. Thus the hypothesis and the conclusion concern two different dimensions.
\end{remark}

\begin{proposition}\label{prop:compression}
If neither $h$ nor $h^+$ has a unique factorization, then there is a homogeneous binomial $f\in I_H$, belonging to a minimal homogeneous generating set, whose class spans both $V_h$ and $V_{h^+}$ as $k$-vector spaces; that is,
\[
V_h=V_{h^+}=k\overline f.
\]
Consequently, if $d=\deg f$, then $d\le_Hh$ and $d\le_Hh^+$, and there are exactly two monomials in $S$ of degree $d$; these monomials are relatively prime.
\end{proposition}

\begin{proof}
By Lemma~\ref{lem:divisor-algebras}, both $V_h$ and $V_{h^+}$ are nonzero. Proposition~\ref{prop:boundary-products} shows that $A_1A_2$ has nonzero components in degrees $\sigma+h+\alpha$ and $\sigma+h+3\alpha$. Since these degrees are distinct,
$\dim_k(A_1A_2)\ge2$.
Let $\m=(t^a\mid a\in H\setminus\{0\})$ be the homogeneous maximal ideal of $R$. Localizing $S$ at $\n$ and $R$ at $\m$ gives an isomorphism
\[
\Tor^{S_{\n}}(R_{\m},k)\cong\Tor^S(R,k)=A
\]
that respects multiplication. Since $S_{\n}$ is regular of dimension four and $R_{\m}$ is Cohen--Macaulay of dimension one, the Auslander--Buchsbaum formula gives $\operatorname{pd}_{S_{\n}}R_{\m}=3$. Put
\[
K=\Ann_{A_1}(A_2).
\]
Proposition~\ref{prop:pd-three} gives $\dim_k(A_1/K)=1$. Moreover, $K$ is graded. To see this, let $\xi=\sum_d\xi_d\in K$ be the decomposition into homogeneous components and let $\eta\in A_2$ be homogeneous. Then $0=\xi\eta=\sum_d\xi_d\eta$. If $d\ne d'$, then $\xi_d\eta$ and $\xi_{d'}\eta$ have different degrees. Hence $\xi_d\eta=0$ for every $d$, and every $\xi_d$ belongs to $K$.

Set $U=V_h+V_{h^+}$. For each $d\in\mathbb Z$,
\[
[U]_d=
\begin{cases}
[A_1]_d,&d\le_Hh\text{ or }d\le_Hh^+,\\
0,&\text{otherwise}.
\end{cases}
\]
If $[U]_d\ne0$, Proposition~\ref{prop:boundary-products} shows that no nonzero element of $[U]_d$ annihilates $A_2$. Hence $[U]_d\cap[K]_d=0$. Since $U$ and $K$ are graded, $U\cap K=0$, and the quotient map restricts to an injection
$U\hookrightarrow A_1/K$.
Thus $\dim_kU\le1$. Since $V_h$ and $V_{h^+}$ are nonzero, we have
\[
V_h=V_{h^+}=U,
\qquad \dim_kU=1.
\]
Choose a homogeneous binomial $f$ in a minimal homogeneous generating set of $I_H$ such that $\overline f$ spans $U$. Lemma~\ref{lem:divisor-algebras} shows that $d=\deg f$ satisfies $d\le_Hh$ and $d\le_Hh^+$ and that there are exactly two monomials in $S$ of degree $d$; these monomials are relatively prime.
\end{proof}

We also record a consequence of the type bound for three-generated numerical semigroups.

\begin{lemma}\label{lem:three-generated-obstruction}
Let $H'=\langle b_1,b_2,b_3\rangle$ be a numerical semigroup, where the $b_i$ are positive integers and need not form a minimal generating set. Let $y,\beta\in\mathbb Z$ with $\beta>0$. Assume that $\beta,2\beta\notin H'$ and that $y,y+\beta,y+2\beta\notin H'$. For any two distinct indices $i,j\in\{1,2,3\}$, there exist $r \in \{0,1,2\}$ and $\ell \in \{i,j\}$ such that $y + b_\ell + r\beta \notin H'$.
\end{lemma}

\begin{proof}
By symmetry, take $i=1$ and $j=2$. Suppose that
\[
y+b_\ell+r\beta\in H'
\qquad(\ell=1,2,\ 0\le r\le2).
\]
For $r=0,1,2$, choose the largest $q_r\ge0$ such that
$z_r=y+r\beta+q_rb_3\notin H'$.
Then $z_r+b_\ell\in H'$ for $\ell=1,2,3$, so $z_r\in\PF(H')$. If $z_r=z_{r'}$ for $r<r'$, then $(r'-r)\beta$ is a positive multiple of $b_3$, contrary to $\beta,2\beta\notin H'$. Hence $z_0,z_1,z_2$ are distinct elements of $\PF(H')$, and $H'$ has type at least three. A numerical semigroup generated by at most three elements has type at most two by Herzog's structure theorem \cite{Herzog}; see also \cite[Section~2]{NariNumataWatanabe}. This contradiction proves the claim.
\end{proof}

\begin{proposition}\label{prop:pure-power}
If neither $h$ nor $h^+$ has a unique factorization, then, after reordering variables, the binomial in Proposition~\ref{prop:compression} can be chosen in the form
$f=X_1^u-X_2^v$ with $u,v>0$.
\end{proposition}

\begin{proof}
Put $d=\deg f$ and write
$f=m_1-m_2$,
where $m_1$ and $m_2$ are monomials. They are relatively prime and are the only monomials of degree $d$, so no variable divides both $m_1$ and $m_2$. Suppose that at least three variables divide $m_1m_2$. After interchanging $m_1$ and $m_2$ and reordering variables if necessary, we may assume that
$X_1\mid m_1$ and $X_2X_3\mid m_2$.
Therefore
\[
d-a_1,d-a_2,d-a_3\in H,\qquad
d-a_1-a_2,d-a_1-a_3\notin H.
\]
Indeed, if $d-a_1-a_2\in H$, then a monomial of degree
$d-a_1-a_2$, multiplied by $X_1X_2$, would give a third monomial of degree
$d$. The same argument applies to $d-a_1-a_3$.
For $1\le r\le3$, put $y_r=h-d+a_1+r\alpha$. Since $d-a_1\in H$, Lemma~\ref{lem:boundary-identities}~\textup{(2)} gives $y_r\notin H$. Since $d-a_1-a_2,d-a_1-a_3\notin H$, the final assertion of Lemma~\ref{lem:boundary-identities}~\textup{(2)} gives
\[
y_r+a_2,y_r+a_3\in \langle a_2,a_3,a_4\rangle.
\]
Let $g=\gcd(a_2,a_3,a_4)$ and $H'=\langle a_2/g,a_3/g,a_4/g\rangle$. Since $g$ divides $a_2$ and $y_r+a_2$, it divides every $y_r$, and hence $g\mid\alpha$. Put $\beta=\alpha/g$. Then
\[
\frac{y_r}{g}\notin H'\qquad(1\le r\le3),
\]
and the three integers $y_1/g,y_2/g,y_3/g$ form an arithmetic progression with common difference $\beta$. Moreover,
\[
\frac{y_r}{g}+\frac{a_2}{g}\in H',
\qquad
\frac{y_r}{g}+\frac{a_3}{g}\in H'
\qquad(1\le r\le3).
\]
If $\beta\in H'$, then $\alpha=g\beta\in H$, and if $2\beta\in H'$, then $2\alpha=g(2\beta)\in H$. The inclusions $\alpha\in H$ and $2\alpha\in H$ contradict Lemma~\ref{lem:semigroup}~\textup{(1)}. Thus $\beta,2\beta\notin H'$, contrary to Lemma~\ref{lem:three-generated-obstruction}. Therefore at most two variables divide $m_1m_2$. Since $m_1$ and $m_2$ are nonconstant and relatively prime, we may reorder the variables so that $m_1=X_1^u$ and $m_2=X_2^v$ for some $u,v>0$.
\end{proof}

\begin{lemma}\label{lem:binary-ci}
Let $x\in H$, and let $f=X_1^u-X_2^v$, where $u,v>0$, be a homogeneous binomial belonging to a minimal homogeneous generating set of $I_H$. If the class of $f$ spans $V_x$ as a $k$-vector space, that is, if $V_x=k\overline f$, then $J_x$ is a complete intersection.
\end{lemma}
\begin{proof}
Put $d=\deg f=ua_1=va_2$. Lemma~\ref{lem:divisor-algebras} gives $[I_H]_x=[(f)]_x$. Suppose that
\[
x=\sum_{i=1}^4u_i a_i=\sum_{i=1}^4v_i a_i.
\]
The corresponding monomials have the same image in
$S/(f)=\bigl(k[X_1,X_2]/(f)\bigr)[X_3,X_4]$.
Therefore $u_3=v_3$ and $u_4=v_4$.
Thus there are integers $c_3,c_4\ge0$ such that every factorization
$x=\sum_{i=1}^4u_i a_i$ satisfies $u_3=c_3$ and $u_4=c_4$.
Since $0\ne\overline f\in V_x$, we have $x-d\in H$.
Choose a monomial of degree $x-d$. Its product with $X_1^u$ has
degree $x$, so the chosen monomial has the form $MX_3^{c_3}X_4^{c_4}$
with $M$ a monomial in $k[X_1,X_2]$. In particular,
\[
 x=d+\deg M+c_3a_3+c_4a_4 .
\]

We claim that
\[
 J_x=(f,\,MX_1X_2,\,X_3^{c_3+1},\,X_4^{c_4+1})=:\mathfrak b.
\]
A monomial survives in $B_x$ exactly when it divides a monomial of
degree $x$. The fixed exponents therefore put the last two powers in
$J_x$. By Lemma~\ref{lem:divisor-algebras}, the only monomials of
degree $d$ are $X_1^u$ and $X_2^v$. It follows that
$d-a_1-a_2\notin H$, since a monomial of that degree multiplied by
$X_1X_2$ would give a third one. If a monomial of degree $x$ were
divisible by $MX_1X_2$, its quotient by
$MX_1X_2X_3^{c_3}X_4^{c_4}$ would be a monomial of degree
$d-a_1-a_2$. Thus $MX_1X_2\in J_x$ as well, and $f\in I_H\subseteq J_x$.

The polynomials $f$ and $MX_1X_2$ are relatively prime in $k[X_1,X_2]$
and form a regular sequence there. Over $k[X_1,X_2]/(f,MX_1X_2)$,
the powers of the independent variables $X_3$ and $X_4$ also form a
regular sequence. Thus $C=S/\mathfrak b$ is an Artinian complete intersection of socle degree
\[
 d+\deg(MX_1X_2)+(c_3+1)a_3+(c_4+1)a_4-\sigma=x .
\]
The graded surjection $C\to B_x$ is nonzero in degree $x$ since
$[B_x]_x=k\overline{t^x}$. It is therefore injective on the
one-dimensional socle $[C]_x$. Every nonzero ideal in an Artinian
local ring meets the socle, so the kernel is zero. This proves
$J_x=\mathfrak b$ and the assertion.
\end{proof}

\begin{proof}[Proof of Theorem~\ref{thm:intro-four}]
It remains to prove \textup{(AP)} $\Rightarrow$ \textup{(PP)} in Conjecture~\ref{conj:main}. Suppose that neither $h$ nor $h^+$ has a unique factorization. Proposition~\ref{prop:compression} gives a homogeneous binomial $f$ belonging to a minimal homogeneous generating set of $I_H$ such that
$V_h=V_{h^+}=k\overline f$.
By Proposition~\ref{prop:pure-power}, after reordering variables, we may choose this binomial in the form $f=X_1^u-X_2^v$ with $u,v>0$. Lemma~\ref{lem:binary-ci}, applied to $x=h$ and $x=h^+$, shows that $J_h$ and $J_{h^+}$ are complete intersections. Theorem~\ref{thm:boundary-reconstruction} then shows that $h$ and $h^+$ have unique factorizations, a contradiction. Hence at least one of $h$ and $h^+$ has a unique factorization. Corollary~\ref{cor:endpoints-uf} now gives \textup{(PP)}. The implications \textup{(PP)} $\Rightarrow$ \textup{(MP)} and \textup{(MP)} $\Rightarrow$ \textup{(AP)} hold in every embedding dimension by Section~\ref{sec:prelim}.
\end{proof}

\begin{remark}\label{rem:four-barriers}
The identities in Lemma~\ref{lem:boundary-identities}, Poincar\'e duality for the Artinian Gorenstein quotients of $R$ associated to $h$ and $h+n\alpha$, the Tor-vanishing criterion, and Theorem~\ref{thm:boundary-reconstruction} work in arbitrary embedding dimension. The proof of Theorem~\ref{thm:intro-four} uses $n=4$ in two places. First, in the contradiction argument, the failure of unique
factorization at both endpoints $h$ and $h^+$ gives $\dim_k(A_1A_2)\ge 2$ by
Proposition~\ref{prop:boundary-products}. Since $R_{\m}$ has projective dimension three over $S_{\n}$, Proposition~\ref{prop:pd-three}
then gives $\dim_k(A_1/\Ann_{A_1}(A_2))=1$.
Second, after one generator is deleted, dividing the remaining three
by their gcd gives a numerical semigroup of type at most two.
\end{remark}

\bibliographystyle{amsplain-arxiv}
\bibliography{refs}
\end{document}